\documentclass[11pt,leqno]{article}

      \usepackage{latexsym}
      \usepackage[centertags]{amsmath}
\usepackage{amsfonts}
\usepackage{amssymb}

\makeatletter\@addtoreset{equation}{section} \makeatother

\renewcommand\thefigure{\thesection.\@arabic\c@figure}
\renewcommand\thetable{\thesection.\@arabic\c@table}

      \newtheorem{theorem}{Theorem}[section]

      \newtheorem{Prop}[theorem]{Proposition}
      
      \def\ni{\noindent}
      
      \def\rf#1{\mbox{$(\ref{#1})$}}

      \def\be{\begin{equation}} 
      \def\ee{\end{equation}} 
      \def\beqn{\begin{eqnarray}} 
      \def\eeqn{\end{eqnarray}} 
      \def\beq{\begin{eqnarray*}} 
      \def\eeq{\end{eqnarray*}}
      \def\proof{{\noindent\bf Proof\quad}\ }
      \def\mb{\mbox} 
      \def\ra{\rightarrow} 

\usepackage{amsmath, amsfonts, amssymb}
\usepackage{mathrsfs}
\usepackage{latexsym}

\newcommand{\scr}[1]{\mathscr #1}

\def\R{\mathbb R}  \def\ff{\frac} 
\def\N{\mathbb N}   
\def\dd{\delta} \def\DD{\Delta}  \def\rr{\rho}
\def\<{\langle} \def\>{\rangle} \def\GG{\Gamma} \def\gg{\gamma}
   \def\pp{\partial}
\def\d{\text{\rm{d}}} \def\bb{\beta} \def\aa{\alpha} \def\D{\scr D}
\def\E{\scr E}   \def\beq{\begin{equation}}
\def\beg{\begin} \def\beq{\begin{equation}}  \def\F{\scr F}
 
\def\e{\text{\rm{e}}}   
     \def\tt{\tilde} 
  
  \def\ee{\varepsilon}

\def\L{\scr L}

      \def\be{\begin{equation}} 
      \def\ee{\end{equation}} 
      \def\beqn{\begin{eqnarray}} 
      \def\eeqn{\end{eqnarray}} 
      \def\beq{\begin{eqnarray*}} 
      \def\eeq{\end{eqnarray*}}
      \def\proof{{\bf Proof:}\ }
      \def\mb{\mbox} 
      \def\ra{\rightarrow} 
      
      \def\rf#1{\mbox{$(\ref{#1})$}}

\begin{document}

\title{ A Class of Infinite Dimensional Diffusion Processes with Connection to Population Genetics}

 \author{Shui Feng\\McMaster
      University\\ shuifeng@mcmaster.ca\and Feng-Yu Wang\\Beijing Normal University and Swansea University\\F.Y.Wang@swansea.ac.uk}

\date{}
\maketitle

\begin{abstract}
Starting from a sequence of independent Wright-Fisher diffusion
processes on $[0,1]$, we construct a class of reversible infinite
dimensional diffusion processes on $\DD_\infty:= \{{\bf x}\in
[0,1]^\N:\ \sum_{i\ge 1} x_i=1\}$ with GEM distribution as the
reversible measure. Log-Sobolev inequalities are established for
these diffusions, which lead to the exponential convergence to the
corresponding reversible measures in the entropy. Extensions are
made to a class of measure-valued processes over an abstract space
$S$. This provides a reasonable alternative to the Fleming-Viot
process which does not satisfy the log-Sobolev inequality when $S$
is infinite
 as observed by W. Stannat \cite{S}.
\end{abstract}

\vspace*{.125in} \noindent {\bf Key words:} Poisson-Dirichlet
distribution, GEM distribution, Fleming-Viot process, log-Sobolev
inequality.
      \vspace*{.125in}


\section{Introduction}

      Population genetics is concerned with the distribution and evolution of gene frequencies in a large
      population at a particular locus. The infinitely-many-neutral-alleles model describes the evolution of the
      gene frequencies under generation independent mutation, and resampling. In statistical
      equilibrium the distribution of gene frequencies is the well known Poisson-Dirichlet distribution introduced
      by Kingman \cite{Kingman75}. When a sample of size $n$ genes is selected from a Poisson-Dirichlet population,
      the distribution of the corresponding allelic partition is given explicitly by the {\it Ewens sampling
      formula}. This provides an important tool in testing neutrality of a population.

      Let
      \[
      \DD_{\infty}=\{{\bf x}=(x_1,x_2,... )\in
[0,1]^\N:
      \sum_{k=1}^{\infty}x_k=1\},
      \]
      and
      \[
      \nabla=\{{\bf x}=(x_1,x_2,... )\in
[0,1]^\N: x_1 \geq x_2\geq \cdots \geq 0,
      \sum_{k=1}^{\infty}x_k=1\}.
      \]
      \noindent  The Poisson-Dirichlet distribution with parameter $\theta
      >0$ (henceforth $PD(\theta)$) is a probability measure
      $\Pi_{\theta}$ on $\nabla$.  We use ${\bf
      P}(\theta)=(P_1(\theta), P_2(\theta), ...)$ to denote the
      $\nabla$-valued random variable with distribution $\Pi_{\theta}$. The component $P_k(\theta)$ represents the
      proportion of the $k$th most frequent alleles. If $u$ is the
      individual mutation rate and $N$ is the effective population size, then the
      parameter $\theta =4N u$ is the population mutation rate.  A different way of describing the distribution
      is through the following size-biased sampling. Let $U_k,
      k=1,2,...$, be a sequence of independent, identically distributed
      random variables with common distribution $Beta(1,\theta)$, and
      set
      \be \label{GEM1}
      X^{\theta}_1 = U_1, X^{\theta}_n = (1-U_1)\cdots (1-U_{n-1})U_n, n \geq 2.
      \ee
      Clearly $(X^{\theta}_1, X^{\theta}_2,\ldots)$ is in space
      $\DD_{\infty}$. The law of $X^{\theta}_1, X^{\theta}_2,...$ is called the one
      parameter GEM distribution and is denoted by
      $\Pi_{\theta}^{gem}$. The descending order of $X^{\theta}_1, X^{\theta}_2,...$ has
      distribution $\Pi_{\theta}$. The sequence $X^{\theta}_k, k=1,2,...$ has the same
      distribution as the
      size-biased permutation of $\Pi_{\theta}$.

       Let $\xi_k,
      k=1,...$ be a sequence of i.i.d. random variables with common
      diffusive distribution $\nu$ on $[0,1]$, i.e., $\nu(x)=0$ for
      every $x$ in $[0,1]$. Set
      \be \label{DIRI1}
      \Theta_{\theta, \nu}=\sum_{k=1}^{\infty}P_k(\theta)\delta_{\xi_k}.
      \ee

      It is known that the law of $\Theta_{\theta,\nu}$ is $Dirichlet(\theta,\nu)$
      distribution, and is the reversible distribution of the
      Fleming-Viot process with mutation operator (cf. \cite{EK81})
      \be\label{parent}
      Af(x)=\frac{\theta}{2}\int_{0}^1 (f(y)-f(x))\nu(dx).
      \ee

       For $0\leq \alpha <1, \theta > -\alpha$, let $\{V_k: k=1,2,...\}$ be a sequence of independent random
      variables such that $V_k$ is a $Beta(1-\alpha, \theta + k\alpha )$ random variable for each $k$. Set
      \be\label{GEM2}
      X^{\theta,\alpha}_1 = V_1, X^{\theta,\alpha}_n = (1-V_1)\cdots (1-V_{n-1})V_n, n \geq 1.
      \ee

       The law of $X^{\theta,\alpha}_1, X^{\theta,\alpha}_2,...$ is called the two-parameter
      GEM distribution and is denoted by $\Pi_{\alpha,\theta}^{gem}$.  The law of the descending order
      statistic of $X^{\theta,\alpha}_1, X^{\theta,\alpha}_2,...$ is called the two-parameter
      Poisson-Dirichlet  distribution (henceforth $\Pi_{\alpha,\theta}$) studied thoroughly in Pitman and Yor
      \cite{PitmanYor97} . The sequence $X^{\theta,\alpha}_k, k=1,2,...$ has the same distribution as the
      size-biased permutation of $\Pi_{\alpha,\theta}$.  It is shown in Pitman
      \cite{Pitman96a} that the two-parameter Poisson-Dirichlet distribution is the most general distribution whose
      size-biased permutation has the same distribution as the GEM representation \rf{GEM2}. A two-parameter ``Ewens
      sampling formula" is obtained in \cite{Pitman95}.
      Let $\Theta_{\theta,\alpha,\nu}$ be defined similarly to $\Theta_{\theta,\nu}$ with $X^{\theta}_k$ being replaced by
      $X^{\theta,\alpha}_k$. We call the law of $\Theta_{\theta,\alpha,\nu}$ a $Dirichlet(\theta,\alpha,\nu)$ distribution.

      The Poisson-Dirichlet distribution and its two-parameter generalization have many similar structures including
      the urn construction in \cite{Hoppe84} and \cite{FengHoppe98}, GEM representation, sampling formula, etc.. But
      we have not seen a stochastic dynamic model similar to the infinitely-many-neutral-alleles model and
      the Fleming-Viot process
      developed for the two-parameter
      Poisson-Dirichlet distribution and $Dirichlet(\theta,\alpha,\nu)$ distribution.

      As the first result in this paper, we are able to construct a class of reversible infinite
      dimensional diffusion processes, the GEM processes, so that both $\Pi_{\theta}^{gem}$ and its two-parameter generalization
      $\Pi_{\alpha,\theta}^{gem}$ appear as the reversible
      measures for appropriate parameters.

      In \cite{S}, the log-Sobolev inequality is studied for the Fleming-Viot
      process with motion given by \rf{parent}. It turns out that
      the log-Sobolev inequality holds only when the type space is
      finite. In the second result of this paper, we will first construct
      a measure-valued process that has the $Dirichlet(\theta,\nu)$ distribution as reversible
      measure. Then we will establish the log-Sobolev inequality for the process.

      The rest of the paper is organized as follows. The GEM
      processes associated with $\Pi_{\theta}^{gem}$ and
      $\Pi_{\alpha,\theta}^{gem}$ are introduced in section 2.
      Section 3 includes the proof of
      uniqueness and the log-Sobolev inequality of the
      GEM process. Finally in section 4, the measure-valued
      process is introduced and the corresponding log-Sobolev
      inequality is established.

      \section{GEM Processes}

      For any $i \geq 1$, let $a_i, b_i$ be two strictly positive numbers. We assume that

      \be \label{boundary}
      \inf_{i}b_i \geq \frac{1}{2}.
      \ee

     Let $X_i(t)$ be  the unique strong
      solution of the stochastic differential equation

      \be\label{SDE1}
      dX_i(t)= (a_i -(a_i+b_i)X_i(t))dt +\sqrt{X_i(t)(1-X_i(t))}dB_i(t), X_i(0) \in [0,1],
      \ee
      where $\{B_i(t): i=1,2,...\}$ are independent one dimensional Brownian motions. It is known that the process
      $X_i(t)$ is reversible with reversible measure $\pi_{a_i,b_i}=Beta(2a_i, 2b_i)$. By direct calculation, the
      scale function of $X_i(\cdot)$ is given by

      \[
      s_i(x)= (\frac{1}{4})^{a_i+b_i}\int_{1/2}^x \frac{dy}{y^{2a_i}(1-y)^{2b_i}}.
      \]

      By \rf{boundary}, we have $\lim_{x \ra 1}s_i(x)=+\infty$ for all $i$. Thus
      starting from the interior of $[0,1]$, the process $X_i(t)$ will not hit
      the boundary $1$ with probability one. Let $E= [0,1)^{\N}$. The process
      \[
      {\bf X}(t)=(X_1(t), X_2(t), ...)
      \]
      is then a $ E$-valued Markov process. Consider the map

      \[
      \Phi:E\ra \bar \DD_\infty,\ \  {\bf x}=(x_1,x_2,...)\ra (\varphi_1({\bf x}),\varphi_2({\bf x}),..)
      \]
      with
      \[
      \varphi_1({\bf x})=x_1, \varphi_n({\bf x}) =x_n (1-x_1)\cdots(1-x_{n-1}), n \geq 2.
      \]

      Clearly $\Phi$ is a bijection and the process ${\bf Y}(t)=\Phi({\bf X}(t))$ is thus a Markov process. Let
      $\bar{E}:= [0,1]^\N$ be the closure of $E$, $C(\bar{E})$ denote the set of all continuous function on
      $\bar{E}$, and $C_{cl}^2(\bar{E})$ be the set of functions in $C(\bar{E})$ that have second order
      continuous derivatives, and depend only on a finite number of coordinates. The sets
$C(E)$ and $C_{cl}^2(E)$ will
      be the respective restrictions of  $C(\bar{E})$ and $C_{cl}^2(\bar{E})$ on $E$.  Then the
      generator of process ${\bf X}(t)$ is given by
      \[
      Lf({\bf x})= \sum_{k=1}^{\infty}\big\{x_k(1-x_k)\frac{\partial^2 f }{\partial x_k^2}
       +(a_k-(a_k+b_k)x_k)\frac{\partial f }{\partial x_k}\big\},\ \  f \in C_{cl}^2(E),
      \]
      and can be extended to $C_{cl}^2(\bar{E})$. The sets $B(E)$ and
      $B(\DD_{\infty})$ are bounded
      measurable functions on $E$ and $\DD_{\infty}$, respectively.

      Let ${\bf a}=(a_1,a_2,\ldots,), {\bf b}=(b_1,b_2,\ldots)$, and
      \[
      \mu_{{\bf a},{\bf b}} = \prod_{k=1}^{\infty}\pi_{a_k,b_k}, \ \ \Xi_{{\bf a},{\bf
      b}}
      =\mu_{{\bf a},{\bf b}} \circ \Phi^{-1}.
      \]
      Then we have
      \begin{theorem}\label{t1}
      The processes ${\bf X}(t)$ and ${\bf Y}(t)$ are reversible with respective reversible measures
      $\mu_{{\bf a},{\bf b}}$ and
      $\Xi_{{\bf a},{\bf b}}$.
       \end{theorem}
      \proof The reversibility of ${\bf X}(t)$ follows from the reversibility of each $X_i(t)$. Now for any two
       $f,g$ in $B(\DD_{\infty})$, the two functions $f\circ \Phi, g\circ \Phi$ are in $B(E)$. From the reversibility of
       ${\bf X}(t)$, we have for any $t >0$,
       \beq
        \int_{\DD_{\infty}}f({\bf y})E_{{\bf y}}[g({\bf y}(t))]\Xi_{{\bf a},{\bf b}}(d{\bf y})
       &=& \int_{E} f(\Phi({\bf x}))E_{{\bf x}}[g(\Phi({\bf x}(t)))]\mu_{{\bf a},{\bf b}}(d{\bf x})\\
       & =& \int_{E} g(\Phi({\bf x}))E_{{\bf x}}[f(\Phi({\bf x}(t)))]\mu_{{\bf a},{\bf b}}(d{\bf
       x})\\
       &=&\int_{\DD_{\infty}}g({\bf y})E_{{\bf y}}[f({\bf y}(t))]\Xi_{{\bf a},{\bf b}}(d{\bf y}).
       \eeq
   Hence ${\bf Y}(t)$ is reversible with reversible measure $\Xi_{{\bf a},{\bf b}}$. $\Box$

\vspace{0.6cm}
 \noindent {\bf Remark.} The one parameter GEM
distribution $\Pi_{\theta}^{gem}$ corresponds to $a_i = \frac{1}{2},
b_i=\frac{\theta}{2}$, and the two parameter GEM distribution
$\Pi_{\alpha,\theta}^{gem}$ corresponds to $a_i=\frac{1-\alpha}{2},
b_i=\frac{\theta +i \alpha}{2}$.

\section{Uniqueness and Poincar\'e/Log-Sobolev Inequalities}

Let
$$\bar\DD_\infty:=\{{\bf x}\in [0,1]^\N:\ \sum_{i=1}^\infty x_i\le 1\}$$
be the closure of space $\DD_{\infty}$ in $\R^\N$ under the topology
induced by cylindrically continuous functions. The probability
$\Xi_{{\bf a},{\bf b}}$ can be extended to the space $\bar\DD_\infty$.
For simplicity, the same notation is used to denote this extended
probability measure.

Now, for ${\bf x}\in \bar \DD_\infty$ such that
\[
\sum_{i=1}^n x_i < 1, \ \mb{for all finite}\ n,
\]
let

$$\L ({\bf x})= \sum_{i,j= 1}^\infty a_{ij}({\bf x})\ff{\pp^2}{\pp x_i\pp x_j} +
\sum_{i=1}^\infty b_i({\bf x})\ff{\pp}{\pp x_i},$$ where

\beg{equation*}\beg{split} & a_{ij}({\bf x}):= x_i x_j
\sum_{k=1}^{i\land j} \ff{(\dd_{ki}(1-\sum_{l=1}^{k-1}x_l)-x_k )
(\dd_{kj}(1-\sum_{l=1}^{k-1}x_l)-x_k )}{x_k(1-\sum_{l=1}^k
x_l)},\\
&b_i({\bf x}):= x_i\sum_{k=1}^i
\ff{(\dd_{ik}\big(1-\sum_{l=1}^{k-1}x_l\big)-x_k ) (a_k
\big(1-\sum_{l=1}^{k-1}x_l\big)- (a_k+b_k)x_k )}{x_k(1-\sum_{l=1}^k
x_l)}.\end{split}\end{equation*} Here and in what follows, we set $
\sum_{i=1}^0 =0$ and $\prod_{i=1}^0 =1$ by conventions. By treating
$\frac{0}{0}$ as one, the definition of $\L ({\bf x})$ can be
extended to all points in  $\bar \DD_\infty$. Through direct
calculation one can see that $\L$ is the generator of the GEM
process.

It follows from direct calculation that

\begin{equation}\label{ab} \sum_{i,j=1}^\infty |a_{ij}(x)|\le 3,\ \ \
|b_i(x)|\le \sum_{k=1}^i (b_k x_k+ a_k),\ \ \ x\in \bar \DD_\infty.
\end{equation}
Indeed, since $1-\sum_{l=1}^{i-1} x_l\ge x_i$ and $\sum_{1\le
i<j<\infty} x_ix_j\le \ff 1 2,$  we obtain

\beg{equation*}\beg{split} \sum_{i,j=1}^\infty |a_{ij}(x)|& =
\sum_{i=1}^\infty a_{ii}(x)+  2 \sum_{1\le i<j<\infty}
|a_{ij}(x)|\\
&\le  \sum_{i=1}^\infty x_i^2 \Big( \ff{1-\sum_{l=1}^i x_l }{x_i} +
\sum_{k=1}^{i-1} \ff {x_k}{1-\sum_{l=1}^k x_l}\Big)\\
&\qquad +   2 \sum_{1\le i<j<\infty} x_ix_j \Big(1 +
\sum_{k=1}^{i-1} \ff{x_k}{1-\sum_{l=1}^{k}x_l}\Big)\\
&\le \sum_{i=1}^\infty x_i\Big( 1-\sum_{l=1}^i x_l +
\sum_{k=1}^{i-1}x_k\Big)\\
&\qquad  +   2 \sum_{i=1}^\infty x_i \sum_{j= i+1}^\infty x_j
\Big(1 + \ff{\sum_{k=1}^{i-1} x_k}{\sum_{l=i+1}^\infty x_l}\Big)\\
&\le  1+   2=3.\end{split}\end{equation*} Thus, the first inequality
in (\ref{ab}) holds. Similarly, the second inequality  also holds.

Let $$\GG(f,g)(x)= \sum_{i,j=1}^\infty a_{ij}(x)\ff{\pp f(x)}{\pp
x_i}\ff{\pp g(x)}{\pp x_j}.$$ Then $\GG(f,f)\in
C_b(\bar\DD_{\infty})$ for any  $f\in C^1_b(\bar\DD_{\infty})$.

For each $a>0, b>0$, let $\aa_{a,b}$ be the largest constant such
that for $\ f\in C_b^1([0,1])$ the log-Sobolev inequality

\begin{equation}\label{Log} \pi_{a,b}(f^2\log f^2) \le \ff 1 {\aa_{a,b}} \int_0^1 x
(1-x)f'(x)^2 \pi_{a,b}(\d x) +\pi_{a,b}(f^2)\log \pi_{a,b}(f^2)
\end{equation}
holds. According to \cite[Lemma 2.7]{S}, we have $\aa_{a,b}\ge
\ff{a\land b}{320}.$ Moreover, it is easy to see that for $a,b>0$ the operator

$$r(1-r)\ff{d^2}{d r^2} + (a-(a+b)r)\ff{d}{d r}$$ on $[0,1]$ has a spectral gap $a+b$ with eigenfunction
$h(r):= a-(a+b)r.$ So, the Poincar\'e inequality

\begin{equation}\label{P} \pi_{a,b}(f^2) \le \ff 1 {a+b} \int_0^1 x
(1-x)f'(x)^2 \pi_{a,b}(\d x) +\pi_{a,b}(f)^2\end{equation} holds.

Let $C^{\infty}_{cl}([0,1]^\N)$ denote the set of all bounded,
$C^{\infty}$ cylindrical functions on $[0,1]^\N$, and
$$\scr FC_b^\infty=\{f|_{\bar\DD_{\infty}}: f \in  C^{\infty}_{cl}([0,1]^{\N})\}.$$

Then we have the following theorem.

\beg{theorem}\label{T1.1}  For any $f,g\in \scr FC_b^\infty$, we
have

\begin{equation}\label{A}\E(f,g):=\Xi_{{\bf a},{\bf b}}(\GG(f,g))= -\Xi_{{\bf a},{\bf b}}(f\L g).
\end{equation}
Consequently, $(\E, \scr FC_b^\infty)$ is closable in $L^2(\bar
\DD_\infty; \Xi_{{\bf a},{\bf b}})$ and its closure is a
conservative regular Dirichlet form, which satisfies the Poincar\'e inequality

\begin{equation}\label{1.1'}\Xi_{{\bf a},{\bf b}}(f^2)\le \ff{1}{\inf_{i\ge 1}(a_i+b_i)}\E(f,f),\ \ \ f\in\D(\E), \Xi_{{\bf a},{\bf b}}(f)=0.
\end{equation}
If moreover $\inf\{a_i \wedge b_i: i\geq 1\}
>0,$ the log-Sobolev
inequality

\begin{equation}\label{1.1}\Xi_{{\bf a},{\bf b}}(f^2\log f^2)\le \ff{1}{\beta_{{\bf
a},{\bf b}}}\E(f,f),\ \ \ f\in\D(\E), \Xi_{{\bf a},{\bf b}}(f^2)=1
\end{equation}
holds for some $\beta_{{\bf a},{\bf b}}\ge \inf\{\frac{a_i\wedge b_i
}{320}:i\geq 1\}>0$.
\end{theorem}
\proof For any $f,g\in \scr FC_b^\infty,$ there exists $n\ge 1$ such
that

\begin{equation}\label{c} f({\bf x})= f(x_1,\cdots, x_n),\ g({\bf x})= g(x_1,\cdots, x_n),\
\ {\bf x}=(x_1,\cdots, x_n,\cdots)\in [0,1]^\N.\end{equation} Let
\[
\varphi^{(n)}({\bf x})=(\varphi_1({\bf x}),\ldots,\varphi_n({\bf
x})),
\]
which maps $[0,1]^n$ on to $\DD_n:=\{x\in [0,1]^n:\
\sum_{i=1}^nx_i\le 1\}.$ Define
$$L_n:= \sum_{i=1}^n x_i(1-x_i) \ff{\pp }{\pp x_i^2}
+\sum_{i=1}^n (a_i-(a_i+b_i)x_i)\ff{\pp}{\pp x_i},$$ and
\[
\pi_{{\bf a},{\bf b}}^{n}=\prod_{i=1}^n \pi_{a_i,b_i},\ \ \
\Xi^n=\pi_{{\bf a},{\bf b}}^n \circ {\varphi^{(n)}}^{-1}.
\]

Then, regarding $\{\Xi^n:= \pi_{{\bf a},{\bf
b}}^n\circ{\varphi^{(n)}}^{-1}:\ n\ge 1\}$ as probability measures
on

$\bar\DD_\infty$, by letting $\Xi^n:= \Xi^n(\d x_1\cdots \d
x_{n})\times \dd_0(\d x_{n+1},\cdots)$, it converges weakly to
$\Xi_{{\bf a},{\bf b}}$. Since $L_n$ is symmetric w.r.t. $\pi_{{\bf
a},{\bf b}}^n$ we have

\begin{equation}\label{3.1}\beg{split} &\int_{[0,1]^n} \sum_{i=1}^n x_i(1-x_i)
\Big(\ff{\pp}{\pp x_i} f\circ \varphi^{(n)}\Big)\Big(\ff{\pp}{\pp
x_i} g\circ
\varphi^{(n)}\Big)   \d\pi_{{\bf a},{\bf b}}^n \\
&=-\int_{[0,1]^n} g\circ\varphi^{(n)} L_n f\circ\varphi^{(n)}
\d\pi_{{\bf a},{\bf b}}^n.\end{split}\end{equation}

Noting that

$$\varphi_i({\bf x})= x_i\prod_{l=1}^{i-1}(1-x_l),\ \ x_i=
\ff{\varphi_i({\bf x})}{1-\sum_{l=1}^{i-1} \varphi_l({\bf x})},\ \ \
i\ge 1,$$
 we have

$$\ff{\d f\circ\varphi^{(n)}({\bf x})}{\d x_i} =
\sum_{j\ge i} \ff{(\dd_{ij} - x_i)\varphi_j({\bf
x})}{x_i(1-x_i)}\ff{\d f}{\d \varphi_j}\circ\varphi^{(n)}({\bf
x}).$$ Therefore,

\begin{equation}\label{3.2} \beg{split}&\int_{[0,1]^n} \sum_{i=1}^n x_i(1-x_i)
\Big(\ff{\pp}{\pp x_i} f\circ \varphi^{(n)}\Big)\Big(\ff{\pp}{\pp
x_i} g\circ
\varphi^{(n)}\Big)   \d\pi_{{\bf a},{\bf b}}^n\\
&= \int_{[0,1]^n} \GG (f , g )\circ \varphi^{(n)}\d\pi_{{\bf a},{\bf
b}}^n= \int_{\DD_n} \GG(f,g) \d\Xi^n.\end{split}\end{equation} By
(\ref{ab}) and (\ref{c}), we have $\GG(f,g)\in C_b(\bar
\DD_{\infty})$ so that the weak convergence of $\Xi^n$ to $\Xi_{{\bf
a},{\bf b}}$ implies

\begin{equation}\label{3.3} \lim_{n\to\infty} \int_{\DD_n}\GG(f,
g)\d\Xi^n=\int_{\bar\DD_\infty} \GG(f,g)\d\Xi_{{\bf a},{\bf
b}}.\end{equation} Similarly, by straightforward calculations we
find

$$L_n f\circ\varphi^{(n)}({\bf x}) = (\L f)\circ\varphi^{(n)}({\bf x}).$$
Moreover, (\ref{ab}) and (\ref{c}) imply that  $g\L f\in C_b(\bar
\DD_{\infty}).$  Thus we arrive at

$$\lim_{n\to\infty}\int_{\DD_n} g\circ \varphi^{(n)}L_n f\circ \varphi^{(n)}\d\pi_{{\bf a},{\bf b}}^n=
\int_{\bar\DD_{\infty}} g \L f\d\Xi_{{\bf a},{\bf b}}.$$ Therefore,
(\ref{A}) follows by combining this with (\ref{3.2}) and
(\ref{3.3}). This implies the closability of $(\E, \scr
FC_b^\infty)$, while the regularity of its closure follows from the
compactness of ${\bar\DD_{\infty}}$ under the usual metric

$$\rr(\mathbf x, \mathbf y):= \sum_{i=1}^\infty 2^{-i} |x_i-y_i|.$$
Indeed, it is trivial that $\D(\E)\cap C_0([0,1]^\N)\supset \scr
FC_b^\infty$ which is dense in $\D(\E)$ under  $\E_1^{1/2}$ given by

$$
\E_1(f,f)=\E(f,f)+\|f\|^2_{2}.$$ Moreover,  for any $F\in
C({\bar\DD_{\infty}})= C_0({\bar\DD_{\infty}})$, by its uniform
continuity due to the compactness of the space,   $$\mathbf
{\bar\DD_{\infty}}\ni \mathbf x\mapsto F_n (\mathbf x):=
F(x_1,\cdots, x_n, 0,0,\cdots),\ \ \ n\ge 1$$ is a sequence of
continuous cylindric functions converging uniformly to $F.$ Since a
cylindric continuous function can be uniformly approximated by
functions in $\scr FC_b^\infty$ under the uniform norm, it follows
that $\scr FC_b^\infty$ is dense in $C_0({\bar\DD_{\infty}})$ under
the uniform norm. That is, the Dirichlet form $(\E,\D(\E))$ is
regular.

\

 Next, the desired Poincar\'e and log-Sobolev inequalities can be deduced from (\ref{P}) and (\ref{Log}) respectively.
 For simplicity, we only prove the latter.
By the additivity property of the log-Sobolev inequality (cf.
\cite{G2}),

$$\mu^n (h^2\log h^2)\le \ff 1 {\bb^n_{\mathbf a,\mathbf b}}\int_{[0,1]^n}
\sum_{i=1}^n x_i(1-x_i) \Big(\ff{\pp h}{\pp x_i}\Big)^2\d\pi_{{\bf
a},{\bf b}}^n +\mu^n(h^2)\log \pi_{{\bf a},{\bf b}}^n(h^2)$$ holds
for all $h\in C_b^1([0,1]^n),$ where
\[
\beta^n_{{\bf a},{\bf b}}= \inf\{\alpha_{a_i,b_i}: i=1,\ldots,n\},
f^{(n)}({\bf x})=f(x_1,\ldots,x_n,0,\ldots).
\]
Combining this with (\ref{3.2}), for any $f\in\D$, the domain of
$\L$, we have

$$\Xi^n({f^{(n)}}^2\log {f^{(n)}}^2)\le \ff 1 {\beta^n_{{\bf a},{\bf b}}}\int_{\DD_n}
\GG^{(n)}(f, f) \d\Xi^n + \Xi^n({f^{(n)}}^2)\log
\Xi^n({f^{(n)}}^2).$$
 Therefore, as explained above,
(\ref{1.1}) for $f\in \D$ follows immediately by letting
$n\to\infty.$ Hence, the proof is completed since $\D(\E)$ is the
closure of $\D$ under  $\E_1^{1/2}.$ \hfill $\Box$

\

We remark that since $(\E,\D(\E))$ is regular, according to
\cite{FOT, MR}, $(L,\D)$ generates a Hunt process whose semigroup
$P_t$ is unique in $L^2(\Xi_{{\bf a},{\bf b}}).$ Thus the GEM
process constructed in section 2 is the unique Feller process
generated by $\L$. Moreover, it is well-known that the log-Sobolev
inequality (\ref{1.1}) implies that $P_t$ converges to $\Xi_{{\bf
a},{\bf b}}$ exponentially fast in entropy; more precisely (see e.g.
\cite[Proposition 2.1]{B}),

$$\Xi_{{\bf a},{\bf b}}(P_t f\log P_t f)\le \e^{-4\beta_{{\bf a},{\bf b}} t}\Xi_{{\bf a},{\bf b}}(f\log f),
\ \ \ f\ge 0, \Xi_{{\bf a},{\bf b}}(f)=1.$$ Moreover, due to Gross
\cite{G}, the log-Sobolev inequality is also equivalent to the
hypercontractivity of $P_t$.

Thus, according to Theorem \ref{T1.1}, we have constructed a
diffusion process which converges to its reversible distribution
$\Xi_{{\bf a},{\bf b}}$ in entropy exponentially fast.

\section{Measure-Valued Process}

It was shown in Stannat \cite{S} that the log-Sobolev inequality
fails to hold for the Fleming-Viot process with parent independent
mutation when there are infinite number of types. In this section,
we will construct a class of measure-valued processes for which the
log-Sobolev inequality holds even when the number of types is
infinity.

Let us first consider a measure-valued processes on a Polish space $S$
induced by the above constructed process and a proper Markov process
on $S^\N.$ More precisely, let $X_t:= (X_1(t),\cdots,
X_n(t),\cdots)$ be the Markov process on $\DD_\infty$ associated to
$(\E,\D(\E)$, and  $\xi_t:= (\xi_1(t),\cdots, \xi_n(t),\cdots)$ be a
Markov process on $S^\N$, independent of $X_t$.  We consider the
measure-valued process

$$\eta_t:= \sum_{i=1}^\infty X_i(t)\dd_{\xi_i(t)},$$
where $X_i$ can be viewed as the proportion of the $i$-th family in
the population, and $\xi_i$ its type or label. Then the above
process describes the evolution of all (countably many) families on
the space $S$. Let $\scr M_1$ be the set of all probability measures
on $S$. Then the state space of this process is

$$\scr M_0:= \{\gg\in \scr M_1:\ \ \text{supp}\,\gg\ \text{contains\
at\ most\ countably\ many\ points}\},$$ which is dense in $\scr M_1$
under the weak topology.

Due to Theorem \ref{T1.1}, if $\xi_t$ converges to its unique
invariant probability measure $\nu$ on $S^\N$, then $\eta_t$
converges to $\Pi:= (\Xi_{\mathbf a, \mathbf b}\times\nu)\circ
\psi^{-1}$ for

$$\psi: \DD_\infty\times S^\N\to \scr M_0;\ \ \psi(\mathbf x,\xi):=
\sum_{i=1}^\infty x_i\dd_{\xi_i}.$$ Unfortunately the process
$\eta_t$ is in general non-Markovian. So we like to modify the
construction by using Dirichlet forms.

Let $\nu$ be a probability measure on $S^\N$ and
$(\E_{S^\N},\D(\E_{S^\N}))$ a conservative symmetric Dirichlet form
on $L^2(\nu).$  We then construct the corresponding quadratic form
on $L^2(\scr M_0; \Pi)$ as follows:

\beg{equation*}\beg{split} &\E_{\scr M_0}(F,G):=
\int_{S^\N} \E(F_\xi, G_\xi)\nu(\d \xi) +\int_{\DD_\infty}\E_{S^\N}(F_{\mathbf x}, G_{\mathbf x})\pi_{a,b}(\d {\mathbf x})\\
&F,G\in \D(\E_{\scr M_0}):= \big\{H\in L^2(\Pi):\ H_{\mathbf x}:=
H\circ\psi(x,\cdot)\in \D(\E_{S^\N})\ \text{for}\ \Xi_{{\bf a},{\bf
b}}\text{-a.s.}\
{\mathbf x},\\
&\qquad\qquad\qquad H_\xi:= H\circ \psi(\cdot,\xi)\in \D(\E) \
\text{for}\ \nu\text{-a.s.}\ \xi, \text{such\ that\ } \E_{\scr
M_0}(H,H)<\infty\big\}.\end{split}\end{equation*} Since $\Pi$ has
full mass on $\scr M_0$, to make the state space complete one may
also consider the above defined form a symmetric form on $L^2(\scr
M_1;\Pi) (=L^2(\scr M_0;\Pi))$.

\beg{theorem} \label{T1.2} Assume  there exists $\aa>0$ such that

$$ \nu(f^2\log f^2)\le \ff 1 {\aa} \E_{S^\N}(f,f) +\nu(f^2)\log \nu(f^2),\ \  \ f\in
\D(\E_{S^\N})$$ holds, then

\begin{equation}\label{LogM} \Pi(F^2\log F^2)\le \ff 1 {\aa\land\beta_{{\bf a},{\bf b}}}
\E_{\scr M_0}(F,F) +\Pi(F^2)\log \Pi(F^2),\ \  \ F\in \D(\E_{\scr
M_0}).
\end{equation}
If moreover $\D(\E_{\scr M_0})\subset L^2(\scr M_1;\Pi)$ is dense,
then $(\E_{\scr M_0}, \D(\E_{\scr M_0}))$ is a conservative
Dirichlet form on $L^2(\scr M_0;\Pi)$ so that the associated Markov
semigroup $P_t$ satisfies

\beg{equation}\label{CC'} \Pi(P_t F\log P_t F)\le \Pi(F\log
F)\e^{-(\beta_{\mathbf a, \mathbf b}\land \aa)t},\ \ \ \ t\ge 0,
F\ge 0, \Pi(F)=1,\end{equation} and $(\E_{\scr M_0}, \D(\E_{\scr
M_0}))$ is regular provided so is $(\E_{S^\N}, \D(\E_{S^\N}))$ and
$S$ is compact.\end{theorem} \proof Let

\beg{equation*}\beg{split}\D(\tt \E)= \big\{\tt F\in L^2(\Xi_{{\bf
a},{\bf b}}\times \nu):&\ \tt F(x,\cdot)\in \D(\E_{S^\N})\
\text{for}\ \Xi_{{\bf a},{\bf b}}\text{-a.s.}\
x,\\
& \tt F(\cdot,\xi)\in \D(\E) \ \text{for}\ \nu\text{-a.s.}\ \xi,
\text{such\ that\ } \tt \E (\tt F,\tt F)
<\infty\big\},\end{split}\end{equation*} where

$$\tt \E(\tt F,\tt G):= \int_{\DD_\infty} \E_{S^\N}(\tt F(\mathbf x,\cdot),
\tt G(\mathbf x,\cdot))\Xi_{{\bf a},{\bf b}}(\d \mathbf x)
+\int_{S^\N} \E(\tt F(\cdot, \xi), \tt G(\cdot, \xi))\nu(\d \xi).$$
Then $(\tt\E,\D(\tt\E))$ is a symmetric Dirichlet form on
$L^2(\DD_\infty\times S^\N; \Xi_{{\bf a},{\bf b}}\times \nu)$ and
(see e.g. \cite[Theorem 2.3]{G2})

\beg{equation}\label{**}(\Xi_{{\bf a},{\bf b}}\times \nu)(\tt
F^2\log \tt F^2)\le \ff 1 {\beta_{{\bf a},{\bf b}}\land \aa}
(\Xi_{{\bf a},{\bf b}}\times \nu)(\tt F^2),\ \ \ \tt F\in \D(\tt
\E), (\Xi_{{\bf a},{\bf b}}\times \nu)(\tt F^2)=1.\end{equation} Let
$\tt P_t$ be the Markov semigroup associated to $(\tt\E,\D(\tt
\E))$. Then $(\ref{CC'})$ follows from the fact that
$\eta_t=\psi(X(t), \xi(t))$ and (\ref{**}) implies (cf.
\cite[Proposition 2.1]{B})

$$
(\Xi_{{\bf a}, {\bf b}}\times \nu)(\tt P_t G\log \tt P_t G)  \le
(\Xi_{{\bf a}, {\bf b}}\times \nu)(G\log G) \e^{-4 (\beta_{\mathbf
a, \mathbf b}\land \aa)t}$$ for all $t\ge 0$ and nonnegative
function $G$ with $(\Xi_{{\bf a}, {\bf b}}\times \nu)(G)=1.$ Since
$F\in \D(\E_{\scr M_0})$ if and only if $F\circ \psi\in \D(\tt \E)$,
and

$$\E_{\scr M_0}(F,F) =\tt\E(F\circ \psi, F\circ \psi),$$  (\ref{LogM}) follows from (\ref{**}).
By the same reason and noting that $(\tt\E,\D(\tt E))$ is a
Dirichlet form, we conclude that $(\E_{\scr M_1}, \D(\E_{\scr
M_0}))$ is a   Dirichlet form provided it is densely defined on
$L^2(\scr M_1; \Pi).$ Finally, if  $S$ is compact then so is $\scr
M_1$ (under the weak topology). Thus, as explained in the proof of
Theorem \ref{T1.1}, for regular $(\E_{S^\N}, \D(\E_{S^\N}))$ the set

$$\{f(\<\cdot, g_1\>, \cdots, \<\cdot, g_n\>):\ n\ge 1, f\in
C_b^1(\mathbb R^n), g_i\in C(S), 1\le i\le n\}\subset C_0(\scr
M_0)\cap \D(\E_{\scr M_1})$$ is dense both in $C_0(\scr M_1))(=
C(\scr M_1))$ under the uniform norm and in $\D(\E_{\scr M_1})$
under the Sobolev norm.  $\Box$

\paragraph{Remark.} Obviously, we have a similar assertion for the Poincar\'e inequality: if there exists $\lambda>0$ such that

$$ \nu(f^2)\le \ff 1 {\lambda} \E_{S^\N}(f,f) +\nu(f)^2,\ \  \ f\in
\D(\E_{S^\N})$$ holds, then
$$ \Pi(F^2)\le \ff 1 {\lambda\land\inf_{i\ge 1}(a_i+b_i)}
\E_{\scr M_0}(F,F) +\Pi(F)^2,\ \  \ F\in \D(\E_{\scr
M_0}).$$

\vspace{0.5cm}

 To see that the above theorem applies to a class of
measure-valued processes on $S$, we present below a concrete
condition on $ \E_{S^\N}$ such that assertions in Theorem \ref{T1.2}
apply. In particular, it is the case if $\E_{S^\N}$ is the Dirichlet
form of a particle system without interactions.

\begin{Prop} Let $\nu_i$ be the $i$-th marginal distribution of $\nu$
and for a function $g$ on $S$ let $g^{(i)}(\xi):= g(\xi_i), i\ge 1.$
Assume that

$$\scr S_0:= \Big\{g\in C_0(S):\ g^{(i)}\in
\D(\E_{S^\N}),\ \sup_{i\ge 1}\E_{S^\N}(g^{(i)},
g^{(i)})<\infty\Big\}$$ is dense in $C_0(S).$ Then $(\E_{\scr M_0},
\D(\E_{\scr M_0}))$ is a symmetric Dirichlet form. \end{Prop}

\proof Under the assumption  and the fact that
$C^2_{cl}(\DD_\infty)$ is dense in $L^2(\scr M_0;\Pi),$  the set

$$\scr S:=\big\{f(\<\cdot, g_1\>, \cdots, \<\cdot, g_n\>):\ n\ge 1,
f\in C_b^1(\mathbb R^n), g_i\in\scr S_0, 1\le i\le n\big\}$$ is
dense in $L^2(\scr M_0; \Pi).$ Therefore, by Theorem \ref{T1.2} it
suffices to show that $\scr S\subset \D(\E_{\scr M_0});$ that is,
for $F:=f(\<\cdot, g_1\>, \cdots, \<\cdot, g_n\>)\in \scr S$, one
has $F\circ \psi\in \D(\tt \E).$ Let

$$F_m(\mathbf x)= F\Big(\sum_{i=1}^m x_i g_1(\xi_i),\cdots, \sum_{i=1}^m x_i
g_n(\xi_i)\Big),\ \ \ \mathbf x\in \DD_\infty,\ \ \ m\ge 1.$$ Since
for fixed $\xi\in S^\N$,

$$\pp_{x_i} F\circ \psi(\cdot, \xi)(\mathbf x)= \sum_{k=1}^n \pp_k f
g_k(\xi_i),\ \ \ i\ge 1$$ is uniformly bounded, one has $F_m\in
\D(\E)$ and (\ref{ab}) yields

$$\E(F_m, F_m)\le C$$ for some constant $C>0$ and all $m\ge 1$ and
$\xi\in S^\N.$ Thus, $F\circ \psi(\cdot, \xi)\in \D(\E)$ for each
$\xi\in S^\N$ and

\beg{equation}\label{A1} \sup_{\xi} \E(F\circ\psi(\cdot,\xi),
F\circ\psi(\cdot,\xi)) \le C.\end{equation} On the other hand, since
$g_k\in \scr S_0, 1\le k\le n,$ noting that for any $\mathbf x\in
\DD_\infty$

$$|F\circ \psi(\mathbf x, \xi)-F\circ \psi(\mathbf x, \xi')|^2\le
(\sum_{k=1}^n \|\pp_k f\|_\infty )^2\sum_{i=1}^\infty x_i
|g_k(\xi_i)-g_k(\xi_i')|^2,$$ we conclude in the spirit of
\cite[Proposition I-4.10]{MR} that $F\circ \psi(\mathbf x, \cdot)\in
\D(\E_{S^\N})$ and

$$\E_{S^\N}(F\circ \psi(\mathbf x, \cdot), F\circ \psi(\mathbf x,
\cdot))\le C'$$ for some $C'>0$ independent of $\mathbf x.$
Combining this with (\ref{A1}) we obtain $F\circ \psi\in \D(\tt
\E).$ $\Box$

{\bf Acknowledgement}.The research of S. Feng is supported by NSERC
of Canada. The research of F.Y. Wang is  supported by
     NNSFC(10121101),
RFDP(20040027009) and the 973-Project of P.R. China.


\begin{thebibliography}{99}
\footnotesize




\bibitem{B} {\sc Bakry, D.} (1997). \emph{On Sobolev and logarithmic Sobolev
inequalities for Markov semigroups,} ``New Trends in Stochastic
Analysis'' (Editors: K. D. Elworthy, S. Kusuoka, I. Shigekawa),
Singapore: World Scientific.

\bibitem{EK81}
{\sc Ethier, S.~N. and Kurtz, T.~G.} (1981). The
infinitely-many-neutral-alleles diffusion model. {\em Adv. Appl.
Prob.,{\bf 13},} 429-452.


       \bibitem{FengHoppe98}
      {\sc Feng, S. and Hoppe, F.~M.} (1998). Large deviation principles for some random combinatorial structures in
       population genetics and Brownian motion. {\em
       Ann. Appl. Probab.}, Vol. 8, No. 4, 975--994.

      \bibitem{G} {\sc Gross, L.} (1976).  Logarithmic Sobolev inequalities. {\em Amer. J.
      Math.,}{\bf 97}, 1061--1083.

      \bibitem{G2} {\sc Gross, L.} (1993). \emph{Logarithmic Sobolev inequalities
       and contractivity properties of semigroups}, Lecture Notes in Math.
       1563, Springer-Verlag.

       \bibitem{FOT} {\sc Fukushima, M., Oshima, Y. and Takeda, M.}
       (1994).
        \emph{Dirichlet Forms and Symmetric Markov Processes,} Walter de
        Gruyter.

      \bibitem{Hoppe84}
      {\sc Hoppe, F.~M.} (1984). P\'olya-like urns and the Ewens sampling formula. {\em
      Journal of Mathematical Biology}, 20, 91--94.


      \bibitem{Kingman75}
      {\sc Kingman, J.~C.~F.} (1975). Random discrete distributions. {\em J. Roy.
      Statist. Soc. B, \bf 37}, 1-22.



      \bibitem{MR} {\sc Ma, Z.~M. and R\"{o}ckner, M.} (1992). \emph{An introduction to the
      theory of {\rm (}non-symmetric{\rm )} Dirichlet forms,} Berlin: Springer.

      \bibitem{Pitman96a}
      {\sc Pitman, J.} (1996). Random discrete distributions invariant under
      size-biased permutation. {\em Adv. Appl. Probab., \bf 28},
      525-539.

      \bibitem{Pitman95}
      {\sc Pitman, J.} (1995). Exchangeable and partially exchangeable random partitions.
      \newblock {\it {Prob. Theory Rel. Fields}}, {\bf 102}, 145--158.

      \bibitem{PitmanYor97}
      {\sc Pitman, J. and Yor, M.} (1997). The two-parameter Poisson-Dirichlet
      distribution derived from a stable subordinator. {\em Ann.
      Probab.},Vol. 25, No. 2, 855-900.



     \bibitem{S} {\sc Stannat, W.} (2000), On validity of the log-Sobolev
      inequality for symmetric Fleming-Viot operators. {\em Ann.
      Probab.,}
      Vol 28, No. 2, 667--684.

\end{thebibliography}
\end{document}